\input amstex
\magnification=\magstep1 
\baselineskip=13pt
\documentstyle{amsppt}
\vsize=8.7truein \CenteredTagsOnSplits \NoRunningHeads

\def\EE{\bold {E\thinspace}}
\def\Pr{\bold {P\thinspace}}
\def\SS{\Cal S}

\topmatter
 
 \title Smoothed counting of 0-1 points in polyhedra \endtitle 
\author Alexander Barvinok  \endauthor
\address Department of Mathematics, University of Michigan, Ann Arbor,
MI 48109-1043, USA \endaddress
\email barvinok$\@$umich.edu \endemail
\date July 13, 2021 \enddate
\thanks  This research was partially supported by NSF Grant DMS 1855428. 
\endthanks 
\abstract Given a system of linear equations $\ell_i(x)=\beta_i$ in an $n$-vector $x$ of 0-1 variables, we compute 
the expectation of $\exp\left\{- \sum_i \gamma_i \left(\ell_i(x) - \beta_i\right)^2\right\}$, where $x$ is a vector of independent Bernoulli random variables and $\gamma_i >0$ are constants. The algorithm runs 
in quasi-polynomial $n^{O(\ln n)}$ time under some sparseness condition on the matrix of the system. The result is based on the absence of the zeros of the analytic continuation of the expectation for complex probabilities, which can also be interpreted as the absence of a phase transition in the Ising model with a sufficiently strong external field. We discuss applications to (perfect) matchings in hypergraphs and randomized rounding in discrete optimization.
\endabstract
\endtopmatter
\document

\head 1. Introduction and examples \endhead

\subhead (1.1) Linear equations in 0-1 vectors \endsubhead Let $A=\left(\alpha_{ij}\right)$ be an $m \times n$ real matrix and let $b=\left(\beta_1, \ldots, \beta_m\right)$ 
be a real $m$-vector, where $m \leq n$. As is well-known, the problem of finding if there is a solution $\xi_1, \ldots, \xi_n$ to the system of linear equations in 0-1 variables
$$\split \sum_{j=1}^n \alpha_{ij} &\xi_j = \beta_i \quad \text{for} \quad i=1, \ldots, m \\
&\xi_j \in \{0, 1\} \quad \text{for} \quad j=1, \ldots, n \endsplit \tag1.1.1$$
 is NP-hard, while counting all such solutions in a $\#$P-hard problem. Motivated by the general difficulty of the problem and inspired by ideas from statistical physics,
we suggest a way of ``smoothed counting", which in some non-trivial cases turns out to be computationally feasible, at least in theory, and gives some information about 
``near-solutions" that satisfy the equations within a certain error. It also allows us to sharpen the procedure of ``randomized rounding" in discrete optimization. As a by-product, 
we prove the absence of phase transition in the Lee - Yang sense for the Ising model with a sufficiently strong external field.

Let us fix some $\gamma_i > 0$, $i=1, \ldots, m$, interpreted as ``weights" of the equations in (1.1.1).  Suppose further, that $\xi_1, \ldots, \xi_n$ are independent 
Bernoulli random variables, so that 
$$\Pr\left(\xi_j=1\right)=p_j  \quad \text{and} \quad \Pr\left(\xi_j=0\right)=1-p_j \quad \text{for} \quad j=1, \ldots, n \tag1.1.2$$
for some $0 < p_j < 1$. Our goal is to compute the expectation 
$$\EE \exp\left\{ -\sum_{i=1}^m \gamma_i \left(-\beta_i + \sum_{j=1}^n \alpha_{ij} \xi_j \right)^2 \right\}. \tag1.1.3$$
Hence every solution to (1.1.1) is accounted for in (1.1.3) with weight 1, while any other 0-1 vector $\left(\xi_1, \ldots, \xi_n\right)$ is accounted for with a weight that is exponentially small in the 
number of violated constraints and the ``severity" of violation.

Clearly, (1.1.3) is always an upper bound on the probability that $(\xi_1, \ldots, \xi_n)$ is a solution to (1.1.1), and that for larger $\gamma_i$ we get sharper bounds.
Generally, we cannot expect to be able to compute (1.1.3) efficiently for $\gamma_i$ that are too large, since that would lead to an efficient algorithm in a $\#$P-hard problem of counting 0-1 solutions to a system of linear equations. How large $\gamma_i$ we can choose will depend on the sparsity of the system (1.1.1) as well as on the choice of probabilities 
$p_1, \ldots, p_n$. 
The choice of probabilities is motivated by the specifics of the problem. For example, if we pick 
$p_j=k/n$ for all $j$ then the probability distribution concentrates around vectors satisfying $\xi_1 + \ldots + \xi_n=k$, so we zoom in on the solutions of (1.1.1) having approximately 
$k$ coordinates equal 1. We discuss another reasonable choice of probabilities in Section 1.5.

Our main results are stated in Section 2. To make them easier to parse, we demonstrate first some of their corollaries.

\subhead (1.2) Example: perfect matchings in hypergraphs \endsubhead 
Let $H=(V, E)$ be a $k$-hypergraph with set $V$ of vertices and set $E$ of edges. Thus the edges of $H$ are some subsets $s \subset V$ such that $|s| \leq k$. 
A {\it perfect matching} in $H$ is a collection $C \subset E$ of edges $s_1, \ldots, s_m$ such that every vertex $v \in V$ belongs to exactly one edge from $C$. As is well-known, to decide whether $H$ contains a perfect matching is an NP-hard problem if $k \geq 3$ and to count all perfect matchings is a $\#$P-hard problem if $k \geq 2$, cf. Problem SP2 in \cite{A+99} and Chapter 17 of \cite{AB09}. If $k=2$, a fully polynomial randomized approximation scheme was constructed by Jerrum, Sinclair and Vigoda \cite{J+04} in the case of a bipartite graph.

For each edge $s \in E$ we introduce a 0-1 variable $\xi_s$. Then the solutions of the system of equations 
$$\sum_{s:\ v \in s} \xi_s =1 \quad \text{for all} \quad v \in V \tag1.2.1$$
are in one-to-one correspondence with perfect matchings in $H$: given a solution $\left(\xi_s:\ s \in E\right)$ we select those edges $s$ for which $\xi_s=1$. The right hand side of the system is the vector of all 1's, $b=\left(\beta_v: v \in V\right)$, where $\beta_v=1$ for all 
$v \in V$.

Suppose now that the hypergraph $H$ is $k$-{\it uniform}, that is, $|s|=k$ for all $s \in E$ and $\Delta$-{\it regular} for some $\Delta \geq 3$, that is, each vertex $v \in V$ is contained in exactly $\Delta$ edges $s \in E$, which is the case in many 
symmetric hypergraphs, such as Latin squares and cubes, see \cite{LL13}, \cite{LL14}, \cite{Ke18}, \cite{Po18}.  Then $|E|=\Delta |V|/k$ and each perfect matching contains exactly $|V|/k$ edges. Let $C \subset E$ be a random collection 
edges, where each edge $s$ is picked into $C$ independently at random with probability 
$$p_s={1 \over \Delta} \quad \text{for all} \quad s \in E, \tag1.2.2$$
so that the expected number of selected edges is exactly $|V|/k$. For a collection $C \subset E$ of edges and a vertex $v \in V$, let $\#(C, v)$ be the number of edges from $C$ that contain $v$. We pick $\gamma_v = \gamma$ for some $\gamma >0$ and all $v \in V$. Then (1.1.3) can be written as 
$$\EE \exp\left\{ -\gamma \sum_{v \in V}  \left(\#(C, v) -1 \right)^2\right\}. \tag1.2.3$$
In Section 2.4, we show that we can choose 
$$\gamma = {\gamma_0 \over k} \tag1.2.4$$
for some absolute constant $\gamma_0 >0$ so that (1.2.3) can be computed within relative error $0 < \epsilon < 1$ in quasi-polynomial time $|E|^{O(\ln |E| -\ln \epsilon)}$. We show that one can choose $\gamma_0=0.025$ and, if $\Delta$ is large enough, one can choose $\gamma_0=0.17$.

The dependence of $\gamma$ on $k$ in (1.2.4) is likely to be optimal, or close to optimal. Indeed, if we could have allowed, for example, 
$\gamma=\gamma_0/k^{1-\epsilon}$ for some fixed $\epsilon >0$ in (1.2.3), we would have been able to approximate (1.2.3) efficiently with any $\gamma>0$, and hence 
compute the probability of selecting a perfect matching with an arbitrary precision. The bootstrapping is accomplished as follows. Given a hypergraph $H=(V, E)$ and an integer $m>1$, let us construct the hypergraph 
$H_m=(V_m, E_m)$. We have $|V_m|=m|V|$ and the vertices of $V_m$ are the ``clones" of the vertices of $V$, so that each vertex of $V$ has $m$ clones in $V_m$.
Each edge $s \in E$ corresponds to a unique edge $s' \in E_m$ such that $|s'|=m|s|$ and $s'$ consists of the clones of each vertex in $s$. We assign the probabilities $p_{s'}=p_s$.
Thus if $H$ is a $k$-uniform hypergraph 
then $H_m$ is $km$-uniform, and if $H$ is $\Delta$-regular then $H_m$ is also $\Delta$-regular. On the other hand, for a collection $C \subset E$ of edges of $H$ and 
the corresponding collection $C' \subset E_m$, we have 
$$\sum_{v \in V_m} \bigl(\#(C', v)-1\bigr)^2= m \sum_{v \in V}\bigl(\#(C, v)-1\bigr)^2.$$
Hence if we could choose $\gamma=\gamma_0/k^{1-\epsilon}$ in (1.2.3), by applying our algorithm to the hypergraph $H_m$ instead of $H$, 
we would have computed (1.2.3) for $H$ with $\gamma=m^{\epsilon} \gamma_0/k^{1-\epsilon}$, and we could have achieved an arbitrarily large $\gamma$ by choosing $m$ large enough.

The standard method of randomized rounding consists of choosing a random collection $C$ of edges from the probability distribution (1.2.2) in lieu of an ``approximate perfect matching", see, for example, Chapter 5 of \cite{MR95} and Section 4.7 of \cite{AS00}.
For a collection $C \subset V$ of edges, we define the penalty function 
$$f(C)=\sum_{v \in V} \bigl(\#(C, v)-1\bigr)^2,$$
which measures how far $C$ is from a perfect matching.
In Section 2.4 we show that for any given $0 < \epsilon < 1$, one can compute (again, in quasi-polynomial time) a particular collection $C_0 \subset E$ of edges such that 
$$\exp\left\{ - \gamma f(C_0) \right\}\  \geq \ (1-\epsilon) \EE  \exp\left\{ -\gamma   f(C) \right\}.$$
It follows that 
$$\Pr\Bigl\{C: \ f(C) \ \leq \ f(C_0) - \rho \Bigr\} \ \leq \ {e^{-\gamma \rho} \over 1-\epsilon} \quad \text{for any} \quad \rho > 0.$$
For example, the probability that a random collection $C$ outperforms $C_0$ with respect to $f$ by $\delta |V|$ for some $\delta >0$ is exponentially small in $\delta |V|/k$.
We note that for a fixed $\Delta$ and $k$, the function $f(C)$ is a random variable with expectation and variance roughly linear in $|V|$. If we assume that 
$f$ has a roughly Gaussian tail, that is,
$$\Pr\bigl\{C: \ f(C) \ \leq \ a - \delta |V| \bigr\} \ \sim \ e^{ - \kappa \delta^2 |V|} $$
for the median $a$ (roughly linear in $|V|$) sufficiently small (constant) $\delta > 0$ and $\kappa>0$  (which is not unreasonable since $f$ is a sum of weakly dependent random variables), then with high probability $f(C_0)$ is smaller than $f(C)$ by a linear in $|V|$ term.

\subhead (1.3) Example: matchings in hypergraphs \endsubhead Let $H=(V, E)$ be a $k$-uniform $\Delta$-regular hypergraph as in Section 1.2. We are still interested in computing 
(1.2.3), only this time we select each edge $s \in E$ into $C$ with a smaller probability 
$$p_s={\omega \over \Delta} \quad \text{for all} \quad s \in E,$$
for some fixed $0 < \omega < 1$. This time the expected cardinality of $C$ is $\omega |V|/k$, so typically $C$ will not cover all vertices of $H$. We also note that once $|C|$ is fixed, the largest weight 
$$\exp\left\{ - \gamma \sum_{v \in V} \left(\#(C, v)-1\right)^2 \right\}$$
is attained if $C$ is a {\it matching}, that is, no two edges from $C$ share a common vertex. It turns out that we can choose $\gamma=\gamma(\omega)$ so that 
$\gamma(\omega) \longrightarrow +\infty$ as $\omega \longrightarrow 0$ and (1.2.3) can be approximated within relative error $0 < \epsilon < 1$ in quasi-polynomial 
$|E|^{O(\ln |E| - \ln \epsilon)}$ time.  More precisely, in Section 2.5, we show that if $\omega >0$ is small enough, we can choose 
$$\gamma={1 \over k} \ln {1 \over \omega}. \tag1.3.1$$
While a matching exists trivially in any hypergraph, unless NP=RP, there is no polynomial time approximation scheme for counting all matchings if $k \geq 3$ \cite{S+19}. Polynomial time approximation algorithms for $k=3$ and $\Delta=3$ are obtained in \cite{D+14} (randomized) and \cite{S+19} (deterministic). 

For $k=2$, the problem of counting all matchings in a given graph is $\#$P-hard \cite{Va79}, while there exists a fully polynomial randomized approximation scheme \cite{JS89}. A deterministic polynomial time algorithm is known if the maximum degree is fixed in advance \cite{B+07}, \cite{PR17}.

\subhead (1.4) Connections to the Ising model \endsubhead 
Given a real symmetric $n \times n$ matrix $G=\left(g_{kj}\right)$ with zero diagonal and a real vector 
$\left(f_1, \ldots, f_n\right)$, 
the partition function in the general Ising model can be written as 
$$\sum_{\eta_1, \ldots, \eta_n =\pm 1} \exp\left\{ \sum_{1\leq k < j \leq n} g_{kj} \eta_k \eta_j +\sum_{j=1}^n f_j \eta_j \right\}. \tag1.4.1$$
Here the values of $\eta_j=\pm 1$ are interpreted as spins of the $j$-th particle, the numbers $g_{kj}$ describe the interaction of the $k$-th and $j$-th particle (if 
$g_{kj}>0$, the interaction is ferromagnetic and if $g_{kj}<0$, the interaction is antiferromagnetic), and $f_j$ describe the external field, see Chapter 3 of \cite{FV18}.

We can write the expectation (1.1.3) in the form (1.4.1) via the substitution
$$\xi_j={\eta_j +1 \over 2} \quad \text{for} \quad j=1, \ldots, n.$$
Let 
$$\split &g_{kj}=-{1 \over 2} \sum_{i=1}^m \gamma_i \alpha_{ik} \alpha_{ij} \quad \text{for} \quad j \ne k \quad  \text{and} \\ &f_j={1 \over 2} \ln {p_j \over 1-p_j} -\sum_{i=1}^m \gamma_i \alpha_{ij} \left(-\beta_i + {1 \over 2} \sum_{k=1}^n \alpha_{ik}\right). \endsplit \tag1.4.2$$
Then (1.1.3) is equal to (1.4.1) multiplied by the constant factor
$$\left(\prod_{j=1}^n p_j (1-p_j)\right)^{1/2} \exp\left\{-\sum_{i=1}^m \gamma_i \left(-\beta_i + {1 \over 2} \sum_{j=1}^n \alpha_{ij}\right)^2 -{1 \over 4} \sum_{i=1}^m \sum_{j=1}^n \gamma_i \alpha_{ij}^2\right\}.$$
In \cite{BB21} we prove that for any $0 < \delta < 1$, fixed in advance, the value of (1.4.1) can be approximated within relative error $0 < \epsilon < 1$ in quasi-polynomial 
$n^{O_{\delta}(\ln n - \ln \epsilon)}$ time provided 
$$\sum_{j:\ j \ne k} \left| g_{jk}\right| \ \leq \ 1-\delta \quad \text{for} \quad k=1, \ldots, n, \tag1.4.3 $$
where the implicit constant in the ``$O$" notation depends on $\delta$ only.
Geometrically, the condition (1.4.3) means that the Lipschitz constant of the quadratic form 
$$\sum_{1 \leq k < j \leq n} g_{kj} \eta_k \eta_j $$
on the Boolean cube $\{-1, 1\}^n$ endowed with the $\ell^1$-metric does not exceed $1-\delta$ (the condition is essentially sharp, modulo NP $\ne$ BPP hypothesis). To avoid dealing with exponentially large numbers, we assume that the coefficients $f_j$ in (1.4.1) are given as $e^{f_j}$. Other than that, the complexity does not depend on $f_j$.

The result of \cite{BB21} and the connection (1.4.2) allows us to handle certain sparse systems (1.1.1). Namely, let us fix integers $r_i \geq 1$ for $i=1, \ldots, m$, integer $c \geq 1$ and
suppose that the matrix $A=\left(\alpha_{ij}\right)$ contains at most $r_i$ non-zero entries in the $i$-th row and at most $c$ non-zero entries in each column, while all entries satisfy the inequalities 
$$|\alpha_{ij}| \ \leq \ 1 \quad \text{for all} \quad i, j.$$
Let us choose 
$$\gamma_i ={1 \over cr_i} \quad \text{for} \quad i=1, \ldots, m.$$
Then for the coefficients $g_{kj}$ defined by (1.4.2) we have 
$$\sum_{k: \ k \ne j} |g_{jk}| = {1 \over 2} \sum_{i=1}^m {|\alpha_{ij}|  \over c} \left({1 \over r_i}  \sum_{k:\ k \ne j}  |\alpha_{ik}|  \right) \ \leq \ {1 \over 2c} \sum_{i=1}^m |\alpha_{ij}| 
\ \leq \ {1 \over 2}$$
and hence (1.4.3) is satisfied with $\delta=0.5$. Consequently, the expectation (1.1.3) can be approximated in quasi-polynomial time $m^{O(1)} n^{O(\ln n - \ln \epsilon)}$ within any 
given relative error $0 < \epsilon < 1$.

We note that the system of equations (1.2.1) for perfect matchings in a $k$-uniform hypergraph is not sparse in the above sense when $k$ is fixed but $\Delta$ is allowed to grow, and the bounds of Section 1.2 do not follow from \cite{BB21}. 

Given an $n \times n$ symmetric matrix $G=\left(g_{kj}\right)$ with zero diagonal, let $\lambda =\lambda_G$ be the largest eigenvalue of $G$. Then the matrix $G-\lambda I$ is negative semidefinite, and hence we can represent $G$ in the form (1.4.2), for some $n \times n$ matrix $A=\left(\alpha_{ij}\right)$, where we choose $m=n$ and $\gamma_i=1$ for all $i$. Our results of Section 2 can be interpreted as saying that the partition function (1.4.1) in the Ising model with an arbitrary matrix $G$ of interactions can be efficiently approximated, provided the external field  is sufficiently strong, that is, if the values of $|f_j|$ are sufficiently large. We say more about the connection in Section 2.6, and also relate it to the Lee - Yang phase transition.

\subhead (1.5) The maximum entropy distribution \endsubhead Given the system (1.1.1),  let $Q \subset {\Bbb R}^n$ be the polytope 
$$\split Q=\biggl\{\left(x_1, \ldots, x_n\right):\quad \sum_{j=1}^n \alpha_{ij} &x_j =\beta_i \quad \text{for} \quad i=1, \ldots, m \quad \text{and} \\
0\ \leq\  &x_j \ \leq\ 1 \quad \text{for} \quad j=1, \ldots, n \biggr\}. \endsplit$$
We define the {\it entropy function} 
$$H(x) = \sum_{j=1}^n x_j \ln {1 \over x_j} +(1-x_j) \ln {1 \over 1-x_j} \quad \text{where} \quad x=\left(x_1, \ldots, x_n\right)$$
and $0 \leq x_j \leq 1$ for $j=1, \ldots, n$, with the standard agreement that at $x_j=0$ or $x_j=1$ the corresponding terms are $0$.

Suppose that the polytope $Q$ has a non-empty relative interior, that is, contains a point $x=\left(x_1, \ldots, x_n\right)$ where $0 < x_j < 1$ for $j=1, \ldots, n$.

One reasonable choice for the probabilities $p_j$ in (1.1.2) is the maximum entropy distribution obtained as the 
solution $p=x$ to the optimization problem:
$$\text{maximize}\quad H(x) \quad \text{subject to} \quad x \in Q. \tag1.5.1$$
This is a convex optimization problem, for which efficient algorithms are available \cite{NN94}. Let $X=\left(\xi_1, \ldots, \xi_n\right)$ be a vector of independent Bernoulli random variables defined by (1.1.2), where $p=\left(p_1, \ldots, p_n\right)$ is the optimal solution in (1.5.1). 
Then $\EE X \in Q$. Moreover, for every point $x \in \{0, 1\}^n \cap Q$, we have 
$$\Pr(X=x)=e^{-H(p)}$$
and hence we get a bound on the number of 0-1 points in $Q$:
$$\left| \{0, 1\}^n \cap Q \right| = e^{H(p)} \Pr(X \in Q) \ \leq \ e^{H(p)},$$
see \cite{BH10} for details. This bound turns out to be of interest in some situations, see, for example, \cite{PP20}.

In this case, our ``smoothed counting" provides an improvement 
$$\left| \{0, 1\}^n \cap Q \right| \ \leq \ e^{H(p)} \EE \exp\left\{ - \sum_{i=1}^m \gamma_i \left(-\beta_i + \sum_{j=1}^n \alpha_{ij} \xi_j \right)^2 \right\},$$
by frequently an exponential in $m$ factor. For example, the distribution (1.2.2) for $k$-uniform $\Delta$-regular hypergraphs is clearly the maximum entropy distribution, and for fixed $k$ and $\Delta$ we get an $e^{\Omega(|V|)}$ factor improvement, compared to the maximum entropy bound.

\head 2. Methods and results \endhead 

\subhead (2.1) The interpolation method \endsubhead
Given an $m \times n$ matrix $A=\left(\alpha_{ij}\right)$, $m$-vector $b=\left(\beta_i\right)$ and weights $\gamma_i$ as in Section 1.1, we consider the polynomial 
$$P_{A, b, \gamma}(\bold z)=\sum_{\xi_1, \ldots, \xi_n \in \{0, 1\}} z_1^{\xi_1} \ldots z_n^{\xi_n}
 \exp\left\{ -\sum_{i=1}^m \gamma_i \left(-\beta_i + \sum_{j=1}^n \alpha_{ij} \xi_j \right)^2 \right\} \tag2.1.1$$
 in $n$ complex variables ${\bold z}=\left(z_1, \ldots, z_n\right)$, where we agree that $z_j^0=1$.
 Hence the expected value (1.1.3) is written as 
$$\left(\prod_{j=1}^n \left(1-p_j\right) \right) P_{A, b, \gamma}\left({p_1 \over 1-p_1}, \ldots, {p_n \over 1-p_n}\right). \tag2.1.2$$
To compute the value of $P_{A,b, \gamma}$ at a particular point $\left(x_1, \ldots, x_n\right)$ we use the interpolation method, see \cite{Ba16} and \cite{PR17} as general references,
as well as recent \cite{Ga20} and \cite{C+21} for connections with other computational approaches, correlation decay and Markov Chain Monte Carlo.
For the interpolation method to work, one should show that that
$$P_{A, b, \gamma}\left(z x_1, \ldots, z x_n\right) \ne 0$$
for all $z$ in some connected open set $U \subset {\Bbb C}$ containing points $0$ and $1$. We establish a sufficient condition for 
$$P_{A, b, \gamma}({\bold z}) \ne 0$$
for all ${\bold z}=\left(z_1, \ldots, z_n\right)$ in a polydisc
$$|z_j| \ < \ \rho_j \quad \text{for} \quad j=1, \ldots, n.$$
We prove the following main result.
\proclaim{(2.2) Theorem} Suppose that the number of non-zero entries in each column of the matrix $A=\left(\alpha_{ij}\right)$ does not exceed $c$ for some integer $c \geq 1$.
 Given real numbers $\rho_j > 0$ for $j=1, \ldots, n$,
we define 
$$\lambda_j =\rho_j \exp\left\{\sum_{i=1}^m \gamma_i \beta_i \alpha_{ij} \right\} \quad \text{for} \quad j=1, \ldots, n.$$
Suppose that 
$$\lambda_j < 1 \quad \text{for} \quad j=1, \ldots, n$$
and that 
$$\sqrt{\gamma_i} \sum_{j=1}^n {\lambda_j \over 1-\lambda_j} \left| \alpha_{ij}\right|  \ \leq \ {1 \over 2\sqrt{c}} \quad \text{for} \quad i=1, \ldots, m.$$
Then 
$$P_{A, b, \gamma} \left({\bold z}\right) \ne 0$$
as long as 
$$|z_j|  \ < \ \rho_j \quad \text{for} \quad j=1, \ldots, n.$$
\endproclaim

Using Theorem 2.2, we obtain an algorithm.

\subhead (2.3) Computing $P_{A, b, \gamma}$ \endsubhead Let us fix a $0 < \delta < 1$ and let $\rho_1, \ldots, \rho_n$ be as in Theorem 2.2. Then for any 
given $x_1, \ldots, x_n$ such that 
$$|x_j | \ \leq \ (1-\delta) \rho_j \quad \text{for} \quad j =1, \ldots, n$$
and any $0 < \epsilon < 1$, the value of 
$$P_{A, b, \gamma}\left(x_1, \ldots, x_n\right)$$
can be approximated within relative error $\epsilon$ in $n^{O_{\delta}(\ln n - \ln \epsilon)}$ time, where the implicit constant in the ``$O$" notation depends on $\delta$ only.
For that, we define a univariate polynomial 
$$g(z)=P_{A, b, \gamma}\left(z x_1, \ldots, z x_n\right) \quad \text{for} \quad z \in {\Bbb C}.$$
Thus $\deg g =n$, we need to approximate $g(1)$ and by Theorem 2.2 we have 
$$g(z) \ne 0 \quad \text{provided} \quad |z| < {1 \over 1-\delta}.$$
As discussed in \cite{Ba16}, Section 2.2, under these conditions, one can approximate in $g(1)$ within relative error $\epsilon$ in $O(n^2)$ time from the values of the derivatives 
$$g^{(k)}(0) \quad \text{for} \quad k \leq O_{\delta}(\ln n - \ln \epsilon),$$
where we agree that $g^{(0)}=g$. From (2.1.1), we have 
$$g(0)=\exp\left\{ - \sum_{i=1}^m \gamma_i \beta_i^2 \right\}$$
while 
$$g^{(k)}(0)=k! \sum\Sb \xi_1, \ldots, \xi_n \in \{0, 1\} \\ \xi_1+ \ldots + \xi_n=k \endSb \exp\left\{ - \sum_{i=1}^m \gamma_i \left(-\beta_i + \sum_{j=1}^n \alpha_{ij} \xi_j \right)^2 \right\}.$$
The direct enumeration of all 0-1 vectors $\left(\xi_1, \ldots, \xi_n\right)$ with $\xi_1+ \ldots + \xi_n=k$ takes $n^{O(k)}$ time and since 
$k=O_{\delta}(\ln n - \ln \epsilon)$, we get the $n^{O_{\delta}(\ln n - \ln \epsilon)}$ complexity of approximating $g(1)$. Here we assume that for any given $\xi_1, \ldots, \xi_n \in \{0, 1\}$, the computation of 
the expression 
$$\exp\left\{ - \sum_{i=1}^m \gamma_i \left(-\beta_i + \sum_{j=1}^n \alpha_{ij} \xi_j \right)^2 \right\}$$
takes unit time. In the bit model of computation, the complexity of the algorithm acquires an additional factor of 
$$\left(mn+\sum_{i,j} |\gamma_i \alpha_{ij}| + \sum_i |\gamma_i \beta_i| \right)^{O(1)}.$$

We now revisit examples of Sections 1.2--1.4 to see how Theorem 2.2 applies there.

\subhead (2.4) Example: perfect matchings in hypergraphs \endsubhead As in Section 1.2, let $H=(V, E)$ be a $k$-uniform $\Delta$-regular hypergraph with $\Delta \geq 3$. 
Let $A=\left(\alpha_{vs}\right)$ be the $|V| \times |E|$ matrix of the system (1.2.1). Hence $\alpha_{vs} \in \{0, 1\}$, every row of $A$ contains $\Delta$ non-zero entries and every column of $A$ contains $k$ non-zero entries, and all non-zero entries are necessarily 1's. Let $b=\left(\beta_v\right)$ be the vector of the right hand sides of (1.2.1). Hence 
$\beta_v=1$ for all $v$. As in Section 1.2, we intend to choose $\gamma_v=\gamma$ for some $\gamma >0$ and all $v \in V$.

Choosing the probabilities $p_s$ as in (1.2.2), in view of (2.1.2), we need to compute 
$$P_{A, b, \gamma}\left({1 \over \Delta-1}, \ldots, {1 \over \Delta-1}\right). \tag2.4.1$$ 
We choose some $0 < \delta < 1$, to be adjusted later, such that for
$$\rho_s=\rho ={1 \over (1-\delta) (\Delta-1)} \quad \text{we have} \quad \rho \ < \ 1$$
(recall that $\Delta \geq 3$). Our goal is to choose $\gamma >0$, the larger the better, such that 
$$P_{A, b, \gamma}\left(z_s: \ s \in E \right) \ne 0$$
provided 
$$|z_s| \ < \rho \quad \text{for all} \quad s \in E.$$
Then we can approximate (2.4.1) by interpolation in quasi-polynomial time, as discussed in Section 2.2. 

We use Theorem 2.2. We need to choose $\gamma >0$ so that for $\lambda_s=\lambda$  we have 
$$\lambda=\rho e^{\gamma k} \ < \ 1 \quad \text{and} \quad \sqrt{\gamma} \Delta {\lambda \over 1-\lambda} \ \leq \ {1 \over 2 \sqrt{k}},$$
that is,
$$\lambda = {e^{\gamma k} \over (1-\delta)(\Delta-1)} \ < \ 1 \quad \text{and} \quad {\lambda \over 1-\lambda} \ \leq \ {1 \over 2 \Delta \sqrt{\gamma k}}.$$
From the second inequality, we get 
$$\lambda \ \leq \ {1 \over 1+2 \Delta \sqrt{\gamma k}} \ < \ 1$$
and hence 
$$e^{\gamma k} \ \leq \ {(1-\delta)(\Delta-1) \over 1+ 2 \Delta \sqrt{\gamma k}}. \tag2.4.2$$
The right hand side of (2.4.2) is an increasing function of $\Delta$, so to find $\gamma=\gamma(k)$ satisfying (2.4.2) for all $\Delta \geq 3$, it suffices to find such $\gamma$ 
satisfying (2.4.2) for $\Delta=3$. Numerical computations show that if we choose a sufficiently small $\delta >0$,  we can choose 
$$\gamma={0.025 \over  k}.$$
If $\Delta$ is large enough, we can choose 
$$\gamma={0.17 \over k}.$$
It turns out that we can compute a particular collection $C_0$ such that 
$$\exp\left\{ - \gamma \sum_{v \in V} \left( \#(C_0, v) -1 \right)^2 \right\} \ \geq \ (1-\epsilon) \EE \exp\left\{ - \gamma \sum_{v \in V} \left( \#(C_0, v) -1 \right)^2 \right\}$$
also in quasi-polynomial time $|E|^{O(\ln |E| - \ln \epsilon)}$. This reduces to computing a sequence of expressions similar to (1.1.3) by the standard application of the method of conditional expectations, see, for example, Chapter 5 of \cite{MR95}. Indeed, the algorithm allows us to compute the conditional expectation, defined by any set of constraints of the type $\xi_j =0$ or $\xi_j=1$. Imposing a condition of this type reduces the computation of 
(1.1.3) to a similar problem, only with fewer variables and possibly different right hand sides $\beta_i$. We note that since the coefficients of the system (1.2.1) are non-negative, when we condition on $\xi_s=0$ or $\xi_s=1$  for a particular edge $s$,  we replace the system with a similar system where the right hand sides $\beta_v$ can only get smaller. Theorem 2.2 then allows us to keep the same value of $\gamma$.
Successively testing for $s \in E$ the conditions $\xi_s=0$ or $\xi_s=1$, and 
choosing each time the one with the larger conditional expectation (which we compute within relative error $\epsilon/|E|^2$), we compute the desired collection $C_0$, while increasing the complexity roughly by a factor of $|E|$.

\subhead (2.5) Example: matchings in hypergraph \endsubhead Here we revisit the example of Section 1.3. This time, we need to compute 
$$P_{A, b, \gamma}\left({\omega \over \Delta-\omega}, \ldots, {\omega \over \Delta-\omega}\right) \quad \text{where} \quad 0 < \omega < 1.$$
Consequently, it suffices to show that 
$$P_{A, b, \gamma}\left(z_s: \ s \in E\right) \ne 0 \quad \text{whenever} \quad |z_s| \ < \ \rho= { \omega \over (1-\delta) \left(\Delta- 1\right)} \quad \text{for all} \quad s \in E$$
and some fixed $0 < \delta < 1$. Using Theorem 2.2, we conclude that we need to choose $\gamma >0$ so that 
$$\lambda = {\omega e^{\gamma k} \over (1-\delta)(\Delta-1)} \ < \ 1 \quad \text{and} \quad {\lambda \over 1-\lambda} \ < \ {1 \over 2\Delta \sqrt{\gamma k}}.$$
It is now clear that if $0 < \omega < 1$ is small enough, we can choose $\gamma$ defined by (1.3.1).

\subhead (2.6) Connections to the Ising model \endsubhead Here we revisit the connection of Section 1.4. Let $G=\left(g_{kj}\right)$ be an $n \times n$ real symmetric matrix with zero diagonal, which we interpret as the matrix of interactions in the Ising model, cf. (1.4.1). Let $\lambda=\lambda_G$ be the largest eigenvalue of $G$. Then the matrix $G-\lambda I$ is negative semidefinite, and hence we can write the entries 
$g_{kj}$ in the form (1.4.2) for some $n \times n$ matrix $A=\left(\alpha_{ij}\right)$ and $\gamma_i=1$ for all $i$.

Suppose that the number of non-zero entries in each column of $A$ does not exceed some $c \geq 1$. For $j=1, \ldots, n$, let us choose $0 < \rho_j < 1$ such that 
$$\sum_{j=1}^n {\rho_j \over 1-\rho_j} |\alpha_{ij}| \ \leq \ {1 \over 2 \sqrt{c}} \quad \text{for} \quad i=1, \ldots, n.$$
Then by Theorem 2.2, we have 
$$\sum_{\xi_1, \ldots, \xi_n \in \{0, 1\}} z_1^{\xi_1} \cdots z_n^{\xi_n} \exp\left\{- \sum_{i=1}^n \left( \sum_{j=1}^n \alpha_{ij}\xi_j \right)^2 \right\} \ne 0$$
as long as $z_1, \ldots, z_n$ are complex numbers such that 
$$|z_j| \ < \ \rho_j \quad \text{for} \quad j=1, \ldots, n.$$
Using (1.4.2), we conclude that 
$$\sum_{\eta_1, \ldots, \eta_n = \pm 1} \exp\left\{ \sum_{1 \leq k < j \leq n} g_{jk} \eta_k \eta_j + \sum_{j=1}^n f_j \eta_j \right\} \ne 0, $$
where $f_1, \ldots, f_n$ are complex numbers with sufficiently small real parts:
$$\split \Re\thinspace f_j \ < \ &{1 \over 2} \ln \rho_j -{1 \over 2} \sum_{i=1}^n \alpha_{ij} \left(\sum_{k=1}^n \alpha_{ik}\right) \\
= \ &{1 \over 2} \ln \rho_j  + \sum_{k=1}^n g_{kj}. \endsplit$$
This can be interpreted as that there is no phase transition in the Lee - Yang sense, see Section 3.7 of \cite{FV18}, provided the external field is strong enough. For comparison, the classical result of Lee and Yang \cite{LY52} establishes that in the ferromagnetic Ising model (that is, when $g_{kj} \geq 0$ for all $k$ and $j$), there is no phase transition as long as the external field is non-zero.

We prove Theorem 2.2 in Section 3. In Section 4, we make some concluding remarks regarding smoothed counting of integer points.

\head 3. Proof of Theorem 2.2 \endhead 

We start with establishing a zero-free region in what may be considered as a Fourier dual functional. The proof of Proposition 3.1 below is somewhat similar to the proof of Theorem 1.1 in \cite{BR19}.
 In what follows, we denote the imaginary unit by $\sqrt{-1}$, so as to use $i$ for indices.
\proclaim{(3.1) Proposition} For $i=1, \ldots, m$ and $j=1,\ldots, n$ let  $\alpha_{ij}$ be real numbers and let $z_j$ be complex numbers.
Suppose that
$$|z_j| \ \leq \ \lambda_j \quad \text{for} \quad j=1, \ldots, n$$
and some $0 < \lambda_j < 1$
and that
$$\left|i: \ \alpha_{ij} \ne 0 \right| \ \leq \ c \quad \text{for} \quad j=1, \ldots, n$$
and some integer $c \geq 1$, so that the $m \times n$ matrix $A=\left(\alpha_{ij}\right)$ has at most $c$ non-zero entries in each column.

If 
$$\sum_{j=1}^n { \lambda_j  \left|\alpha_{ij}\right|  \over 1-\lambda_j} \ \leq \ {1 \over 2\sqrt{c}} \quad \text{for} \quad i=1, \ldots, m$$
Then 
$$\sum_{\sigma_1, \ldots, \sigma_m=\pm 1} \prod_{j=1}^n \left(1+z_j \exp\left\{ \sqrt{-1} \sum _{i=1}^m 
\alpha_{ij} \sigma_i \right\}\right) \ne 0.$$
\endproclaim

Before we embark on the proof of Proposition 3.1, we do some preparations. 

\subhead (3.2) Preliminaries \endsubhead 
Let $\{-1, 1\}^m$ be the discrete cube of all $m$-vectors 
$x=\left(\sigma_1, \ldots, \sigma_m\right)$,
where $\sigma_i =\pm 1$ for $i=1, \ldots, m$. Let $I \subset \{1, \ldots, m\}$ be a set of indices and let us fix some $\epsilon_i \in \{-1, 1\}$ for all $i \in I$. The 
set 
$$F=\Bigl\{(\sigma_1, \ldots, \sigma_m) \in \{-1, 1\}^m:\quad \sigma_i =\epsilon_i \quad \text{for} \quad i \in I \Bigr\}$$ is called a {\it face} of the cube. The indices 
$i \in I$ are {\it fixed} indices of $F$ and indices $i \in \{1, \ldots, m\} \setminus I$ are its {\it free} indices. We define the {\it dimension} by $\dim F=m-|I|$, the cardinality of the set of 
free indices. Thus a face of dimension $k$ consists of $2^k$ points. The cube itself is a face of dimension $m$, while every vertex $x$ is a face of dimension 0.

For a function $f: \{-1, 1\}^m \longrightarrow {\Bbb C}$ and a face $F \subset \{-1, 1\}^m$, we define 
$$\SS(f; F) =\sum_{x \in F} f(x).$$
Suppose that $i$ is a free index of $F$ and let $F^+ \subset F$ and $F^-\subset F$ be the faces of $F$ defined by the constraint $\sigma_i=1$ and $\sigma_i=-1$ respectively. 
Then 
$$\SS(f; F)=\SS(f; F^+) + \SS(f; F^-).$$
Furthermore, if $\SS(f; F^+) \ne 0$ and $\SS(f; F^-) \ne 0$ and if the angle between non-zero complex numbers $\SS(f; F^+)$ and $\SS(f; F^-)$, considered as vectors 
in ${\Bbb R}^2 = {\Bbb C}$, does not exceed $\theta$ for some $0 \leq \theta < \pi$, we have 
$$\left| \SS(f; F) \right| \ \geq \ \left( \cos {\theta \over 2}\right) \left( \left| \SS(f; F^+)\right| + \left| \SS(f; F^-)\right|\right), \tag3.2.1$$
cf. Lemma 3.6.3 of \cite{Ba16}. The inequality (3.2.1) is easily obtained by bounding the length of $\SS(f; F)$
from below by the length of its orthogonal projection onto the bisector of the angle 
between $\SS(f; F^+)$ and $\SS(f; F^-)$.

More generally, suppose that for every face $G \subseteq F$, every free index $i$ of $G$ and the corresponding faces $G^+$ and $G^-$ of $G$, we have that 
$\SS(f; G^+) \ne 0$, $\SS(f; G^-) \ne 0$ and the angle between the two non-zero complex numbers does not exceed $\theta$. Let $I \subset \{1, \ldots,m\}$ be a set of some free indices 
of $F$. For for an assignment $s: I \longrightarrow \{-1, 1\}$ of signs, let $F^s$ be the face of $F$ obtained by fixing the coordinates $\sigma_i$ with $i \in I$ to $s(i)$. Then 
$$\SS(f; F) = \sum_{s:\ I \longrightarrow \{-1, 1\}} \SS(f; F^s)$$
and iterating (3.2.1) we obtain 
$$\left|\SS(f; F)\right| \ \geq \ \left( \cos {\theta \over 2}\right)^{|I|} \sum_{s:\ I \longrightarrow \{-1, 1\}} \left| \SS(f; F^s)\right|. \tag3.2.2$$
Finally, we will use the inequality
$$\left(\cos {\psi \over \sqrt{c}}\right)^c \ \geq \ \cos \psi  \quad \text{for} \quad 0 \leq \psi \leq {\pi \over 2} \quad \text{and} \quad c \geq 1, \tag3.2.3$$
which can be obtained as follows. Since $\tan \psi$ is a convex function on the interval $(0, \pi/2)$, we have 
$$\sqrt{c}\tan {\psi \over \sqrt{c}} \ \leq \ \tan \psi$$
on the interval. Integrating, we obtain 
$$-c \ln \cos {\psi \over \sqrt{c}} \ \leq \ -\ln \cos \psi \quad \text{for} \quad 0 \leq \psi < \pi/2,$$
which is equivalent to (3.2.3).

\subhead (3.3) Proof of Proposition 3.1 \endsubhead
For a given complex vector ${\bold z}=\left(z_1, \ldots, z_n\right)$, satisfying the conditions of the theorem, we consider the function
$\ell(\cdot; {\bold z}):\ \{-1, 1\}^m \longrightarrow {\Bbb C}$ defined by
$$ \ell(x; {\bold z})=\prod_{j=1}^n \left(1+z_j \exp\left\{ \sqrt{-1} \sum _{i=1}^m 
\alpha_{ij} \sigma_i \right\}\right) $$
for $x=\left(\sigma_1, \ldots, \sigma_m\right)$. 
To simplify the notation somewhat, for a face $F \subset \{-1, 1\}^n$,  we denote 
$\SS\bigl(\ell(\cdot; {\bold z}); F\bigr)$ just by $\SS\bigl(\ell({\bold z}); F\bigr)$.

We prove by induction on $d=0, 1, \ldots, m$ the following statement.
\bigskip
\noindent (3.3.1) Let $F \subset \{-1, 1\}^m$ be a face of dimension $d$. Then $\SS\bigl(\ell({\bold z}); F \bigr) \ne 0$. Moreover, if $\dim F >0$ and if $i$ is a free index of 
$F$ then for the faces $F^+, F^- \subset F$ the angle between complex numbers $\SS\bigl(\ell( {\bold z}); F^+\bigr) \ne 0$ and 
$\SS\bigl(\ell( {\bold z}); F^-\bigr) \ne 0$, considered as vectors in ${\Bbb R}^2={\Bbb C}$, does not exceed 
$$\theta = {2 \pi \over 3 \sqrt{c}}.$$
\bigskip
We obtain the desired result when $F=\{-1, 1\}^m$ is the whole cube.
\bigskip
Since $|z_j| < 1$ for $j=1, \ldots, n$, the statement (3.3.1) holds if $\dim F=0$, and hence $F$ is a vertex of the cube. 

Suppose now that (3.3.1) holds for all faces of dimension $d$ and lower. Let $G \subset \{-1, 1\}^m$ be a face of dimension $d$. Since by the induction hypothesis
$\SS\bigl(\ell({\bold z}); G\bigr)\ne 0$ on the polydisc of vectors ${\bold z}=\left(z_1, \ldots, z_n\right) \in {\Bbb C}^n$, satisfying 
$$|z_j| \ \leq \ \lambda_j \quad \text{for} \quad j=1, \ldots, n, \tag3.3.2$$
we can choose a branch of the function 
$${\bold z} \longmapsto \ln \SS\bigl(\ell({\bold z}); G\bigr).$$
For $j=1, \ldots, n$, let us introduce a function
$h_j(\cdot; {\bold z}): \ \{-1, 1\}^m \longrightarrow {\Bbb C}$ defined by 
$$h_j(x; {\bold z})={\exp\left\{\sqrt{-1} \sum_{i=1}^m  \alpha_{ij} \sigma_i \right\} \over 1+z_j \exp\left\{\sqrt{-1} \sum_{i=1}^m  \alpha_{ij} \sigma_i \right\}}$$
for $x=\left(\sigma_1, \ldots, \sigma_m \right)$.
Hence we have
$$\left| h_j(x; {\bold z})\right| \ \leq \ {1 \over 1-\lambda_j} \quad \text{for all} \quad x \in \{-1, 1\}^m \tag3.3.3$$
and
$$\split &{\partial \over \partial z_j} \ell(x; {\bold z})= 
 \exp\left\{ \sqrt{-1} \sum _{i=1}^m 
\alpha_{ij} \sigma_i \right\}\prod_{k:\  k \ne j} \left(1+z_k \exp\left\{ \sqrt{-1} \sum _{i=1}^m 
\alpha_{ik} \sigma_i \right\}\right) \\ &\qquad  \qquad \quad=h_j(x; {\bold z})\ell(x; {\bold z}). \endsplit$$
Therefore,
$${\partial \over \partial z_j} \ln \SS\bigl(\ell({\bold z}); G\bigr) ={ {\partial \over \partial z_j} \SS\bigl(\ell({\bold z}); G) \over \SS\bigl(\ell({\bold z}); G\bigr)}= 
{ \SS\left({\partial \over \partial z_j} \ell({\bold z}); G\right) \over \SS(\ell({\bold z}); G)}=
{\SS(\ell({\bold z}) h_j({\bold z}); G) \over \SS(\ell({\bold z}); G)},$$
where we use $\SS(\ell({\bold z}) h_j({\bold z}); G)$ as a shorthand for $\SS\bigl(\ell(\cdot; {\bold z}) h_j(\cdot; {\bold z}); G\bigr)$.
Our goal is to bound 
$$\left| {\partial \over \partial z_j} \ln \SS\bigl(\ell({\bold z}); G\bigr)\right| = \left| {\SS(\ell({\bold z}) h_j({\bold z}); G) \over \SS(\ell({\bold z}); G)}\right|, \tag3.3.4$$
which will allow us to bound the angle by which $\SS\bigl(\ell({\bold z}); G\bigr)$ rotates as ${\bold z}$ changes inside the polydisc (3.3.2).

If $d=0$ then by (3.3.3), we have 
$$\left| {\partial \over \partial z_j} \ln \SS\bigl(\ell({\bold z}); G\bigr)\right| \ \leq \ {1 \over 1-\lambda_j}.$$
For an index $j=1, \ldots, n$, let 
$$I_j=\Bigl\{i:\ \alpha_{ij} \ne 0 \quad \text{and $i$ is free in $G$} \Bigr\}.$$
Hence $|I_j| \leq c$. Suppose first that $I_j \ne \emptyset$. For an assignment $s: I_j \longrightarrow \{-1, 1\}$ of signs, let $G^s$ be the face of $G$ obtained by 
fixing $\sigma_i=s(i)$ for all $i \in I_j$. Applying the induction hypothesis to $G$ and its faces,  by (3.2.2) we get
$$\left| \SS\bigl(\ell({\bold z}); G\bigr)\right| \ \geq \ \left( \cos {\theta \over 2} \right)^{|I_j|} \sum_{s:\ I_j \longrightarrow \{-1, 1\}}
\left| \SS\bigl(\ell({\bold z}); G^s \bigr)\right|.$$
On the other hand, the function $h_j({\bold z})$ is constant on every face $G^s$, and hence from (3.3.3), we obtain 
$$\split \left| \SS\bigl(\ell({\bold z}) h_j({\bold z}); G\bigr)\right| \ \leq  &\sum_{s:\ I_j \longrightarrow \{-1, 1\}}
\left| \SS\bigl(\ell({\bold z}) h_j({\bold z}); G^s \bigr)\right| \\ \leq \ &{1 \over 1-\lambda_j}  \sum_{s:\ I_j \longrightarrow \{-1, 1\}}
\left| \SS\bigl(\ell({\bold z}); G^s \bigr)\right|. \endsplit$$
Therefore, by (3.3.4), we obtain the bound
$$\left| {\partial \over \partial z_j} \ln \SS\bigl(\ell({\bold z}); G\bigr)\right| \ \leq \ {1 \over (1-\lambda_j) \cos^{|I_j|} (\theta/2)} \ \leq \ 
{1 \over (1-\lambda_j) \cos^{c} (\theta/2)}.  \tag3.3.5$$
If $I_j=\emptyset$ then $h_j({\bold z})$ is constant on $G$ and from (3.3.3) and (3.3.4) we get 
$$\left| {\partial \over \partial z_j} \ln \SS\bigl(\ell({\bold z}); G\bigr)\right| \ \leq \ {1 \over 1-\lambda_j},$$
so (3.3.5) holds as well.

Now we are ready to complete the induction step. Let $F$ be a face of dimension $d+1 > 0$. Let $i$ be a free index of $F$ and let $F^+, F^- \subset F$ be the faces obtained by fixing 
$\sigma_i=1$ and $\sigma_i =-1$ respectively.
Then $\dim F^+ = \dim F^-=d$ and by the induction hypothesis, we have $\SS\bigl(\ell({\bold z}); F^+\bigr) \ne 0$ and $\SS\bigl(\ell( {\bold z}); F^-\bigr) \ne 0$. We need to prove that the 
angle between $\SS\bigl(\ell({\bold z}); F^+\bigr) \ne 0$ and $\SS\bigl(\ell({\bold z}); F^-\bigr) \ne 0$ does not exceed $\theta$.
To this end, we note that 
$$\split &\SS\bigl(\ell({\bold z}); F^+\bigr) = \SS\bigl(\ell( \widehat{{\bold z}}); F^-\bigr) \\
&\text{where} \quad \widehat{z}_j =e^{2 \sqrt{-1} \alpha_{ij}}  z_j \quad \text{for} \quad j=1, \ldots, n. \endsplit $$
Applying (3.3.5) with $G=F^-$, we conclude that the angle between 
$\SS\bigl(\ell({\bold z}); F^+\bigr) \ne 0$ and $\SS\bigl(\ell({\bold z}); F^-\bigr) \ne 0$ does not exceed 
$$\sum_{j=1}^n { \left| z_j - \widehat{z}_j\right| \over (1-\lambda_j)\cos^c (\theta/2)}  \ \leq \ {2 \over \cos^c (\theta/2)} \sum_{j=1}^n {\lambda_j \left| \alpha_{ij}\right| \over 1-\lambda_j}.
$$
Using (3.2.3), we obtain 
$$\cos^c \left({\theta \over 2}\right) = \cos^c {\pi \over 3 \sqrt{c}} \ \geq \ \cos {\pi \over 3} ={1 \over 2},$$
and hence the angle between $\SS\bigl(\ell({\bold z}); F^+\bigr) \ne 0$ and $\SS\bigl(\ell({\bold z}); F^-\bigr) \ne 0$ does not exceed
$$4 \sum_{j=1}^n {\lambda |\alpha_{ij}| \over 1-\lambda_j} \ \leq \ {2 \over \sqrt{c}} \ < \ {2 \pi \over 3 \sqrt{c}} =\theta,$$
which completes the proof.
{\hfill \hfill \hfill} \qed

The following corollary can be considered as a Fourier dual statement to Proposition 3.1.

\proclaim{(3.4) Corollary} For $i=1, \ldots, m$ and $j=1,\ldots, n$, let $0 < \lambda_j < 1$ and  $\alpha_{ij}$ be real numbers and let $z_j$ be complex numbers.
Suppose that
$$\left|i: \ \alpha_{ij} \ne 0 \right| \ \leq \ c \quad \text{for} \quad j=1, \ldots, n$$
and some integer $c \geq 1$ and that
$$\sum_{j=1}^n { \lambda_j  \left|\alpha_{ij}\right|  \over 1-\lambda_j} \ \leq \ {1 \over 2\sqrt{c}} \quad \text{for} \quad i=1, \ldots, m.$$
Then 
$$\sum_{\xi_1, \ldots, \xi_n \in \{0, 1\}} z_1^{\xi_1} \cdots z_n^{\xi_n} \prod_{i=1}^m  \cos\left( \sum_{j=1}^n \alpha_{ij} \xi_j\right)\ne 0$$
provided 
$$|z_j| \ \leq \ \lambda_j \quad \text{for} \quad j=1, \ldots, n.$$
\endproclaim
\demo{Proof}
We have
$$\split  &\prod_{i=1}^m  \cos\left( \sum_{j=1}^n \alpha_{ij} \xi_j\right)\\&\qquad=2^{-m}  \prod_{i=1}^m \left( \exp\left\{ \sqrt{-1} \sum_{j=1}^n \alpha_{ij} \xi_j \right\} 
+ \exp\left\{ - \sqrt{-1} \sum_{j=1}^n \alpha_{ij} \xi_j \right\} \right) \\
&\qquad=2^{-m}\sum_{\sigma_1 \ldots, \sigma_m =\pm 1} \exp\left\{ \sqrt{-1} \sum_{j=1}^n \sum_{i=1}^m \alpha_{ij} \sigma_i \xi_j \right\}. \endsplit$$ 
Consequently,
$$\split &2^m\sum_{\xi_1, \ldots, \xi_n \in \{0, 1\}} z_1^{\xi_1} \cdots z_n^{\xi_n} \prod_{i=1}^m  \cos\left( \sum_{j=1}^n \alpha_{ij} \xi_j\right)\\=
&\sum_{\sigma_1 \ldots, \sigma_m =\pm 1} \sum_{\xi_1, \ldots, \xi_n \in \{0, 1\}} z_1^{\xi_1} \cdots z_n^{\xi_n} 
\exp\left\{ \sqrt{-1} \sum_{j=1}^n \sum_{i=1}^m \alpha_{ij} \sigma_i \xi_j \right\} \\
=&\sum_{\sigma_1 \ldots, \sigma_m =\pm 1} \prod_{j=1}^n \left(1+ z_j \exp\left\{ \sqrt{-1} \sum_{i=1}^m \alpha_{ij} \sigma_i \right\} \right) \ne 0 \endsplit$$
by Proposition 3.1.
{\hfill \hfill \hfill} \qed
\enddemo

Next, we take a limit in Corollary 3.4

\proclaim{(3.5) Corollary} For $i=1, \ldots, m$ and $j=1,\ldots, n$, let $0 < \lambda_j < 1$ and  $\alpha_{ij}$ be real numbers and let $z_j$ be complex numbers.
Suppose that
$$\left|i: \ \alpha_{ij} \ne 0 \right| \ \leq \ c \quad \text{for} \quad j=1, \ldots, n$$
and some integer $c \geq 1$ and that
$$\sum_{j=1}^n { \lambda_j  \left|\alpha_{ij}\right|  \over 1-\lambda_j} \ \leq \ {1 \over 2\sqrt{c}} \quad \text{for} \quad i=1, \ldots, m.$$
Then 
$$\sum_{\xi_1, \ldots, \xi_n \in \{0, 1\}} z_1^{\xi_1} \cdots z_n^{\xi_n} \exp\left\{ -{1 \over 2} \sum_{i=1}^m \left( \sum_{j=1}^n \alpha_{ij} \xi_j\right)^2 \right\} \ne 0$$
provided 
$$|z_j| \  < \ \lambda_j \quad \text{for} \quad j=1, \ldots, n.$$
\endproclaim 
\demo{Proof} Let $A=\left(\alpha_{ij}\right)$ be the $m \times n$ matrix. For an integer $k > 1$, we define the $(km) \times n$ matrix $A^{(k)}=\left(\alpha_{ij}^{(k)}\right)$ as 
follows. First, we divide each row of $A$ by $\sqrt{k}$ and then copy the resulting row $k$ times. Thus every column of $A^{(k)}$ contains at most $ck$ non-zero entries
$$\left| i: \ \alpha_{ij}^{(k)} \ne 0 \right| \ \leq \ ck \quad \text{for} \quad j=1, \ldots, n$$
and we have
$$\sum_{j=1}^n { \lambda_j  \left|\alpha_{ij}^{(k)}\right|  \over 1-\lambda_j} \ \leq \ {1 \over 2\sqrt{ck}} \quad \text{for} \quad i=1, \ldots, km.$$
Applying Corollary 3.4, we conclude that 
$$\sum_{\xi_1, \ldots, \xi_n \in \{0, 1\}} z_1^{\xi_1} \cdots z_n^{\xi_n} \prod_{i=1}^m  \cos^k\left( {1 \over \sqrt{k}} \sum_{j=1}^n \alpha_{ij} \xi_j\right)\ne 0$$
provided 
$$|z_j| \  < \ \lambda_j \quad \text{for} \quad j=1, \ldots, n.$$
Since 
$$\lim_{k \longrightarrow \infty} \cos^k\left( {1 \over \sqrt{k}} \sum_{j=1}^n \alpha_{ij} \xi_j\right) = \exp\left\{ -{1 \over 2} \left( \sum_{j=1}^n \alpha_{ij} \xi_j \right)^2\right\},$$
By Hurwitz' Theorem, see, for example, Chapter 7 of \cite{Kr01}, the polynomial 
$$p\left(z_1, \ldots, z_n\right)=\sum_{\xi_1, \ldots, \xi_n \in \{0, 1\}} z_1^{\xi_1} \cdots z_n^{\xi_n} \exp\left\{ - {1 \over 2} \sum_{i=1}^m \left( \sum_{j=1}^n \alpha_{ij} \xi_j \right)^2\right\}$$
either has no zeros in the domain
$$|z_j| \ < \ \lambda_j \quad \text{for} \quad j=1, \ldots, n \tag3.5.1$$
or is identically zero there. Since 
$$p\left(0, \ldots, 0\right)=1 \ne 0,$$
we conclude that $p\left(z_1, \ldots, z_n\right) \ne 0$ in the domain (3.5.1).
{\hfill \hfill \hfill} \qed 
\enddemo

Next, we deal with non-homogeneous equations.

\proclaim{(3.6) Corollary}  For $i=1, \ldots, m$ and $j=1,\ldots, n$, let $0 < \lambda_j < 1$, $\beta_i$ and $\alpha_{ij}$ be real numbers and let 
$z_j$ be complex numbers. Suppose that
$$\left|i: \ \alpha_{ij} \ne 0 \right| \ \leq \ c \quad \text{for} \quad j=1, \ldots, n$$
and some integer $c \geq 1$ and that
$$\sum_{j=1}^n { \lambda_j  \left|\alpha_{ij}\right|  \over 1-\lambda_j} \ \leq \ {1 \over 2\sqrt{c}} \quad \text{for} \quad i=1, \ldots, m.$$
Then 
$$\sum_{\xi_1, \ldots, \xi_n \in \{0, 1\}} z_1^{\xi_1} \cdots z_n^{\xi_n} \exp\left\{ -{1 \over 2} \sum_{i=1}^m \left( -\beta_i+\sum_{j=1}^n \alpha_{ij} \xi_j\right)^2 \right\} \ne 0$$
provided 
$$|z_j| \  < \ \lambda_j \exp\left\{-\sum_{i=1}^m \beta_i \alpha_{ij} \right\}  \quad \text{for} \quad j=1, \ldots, n.$$
\endproclaim 
\demo{Proof} We have 
$$\left( -\beta_i+\sum_{j=1}^n \alpha_{ij} \xi_j\right)^2=\beta_i^2 - 2\sum_{j=1}^n \alpha_{ij} \beta_i \xi_j + \left(\sum_{j=1}^n \alpha_{ij} \xi_j\right)^2.$$
Denoting 
$$w_j =z_j \exp\left\{ \sum_{i=1}^m \alpha_{ij} \beta_i \right\} \quad \text{for} \quad j=1, \ldots, n,$$
we write 
$$\split &\sum_{\xi_1, \ldots, \xi_n \in \{0, 1\}} z_1^{\xi_1} \cdots z_n^{\xi_n} \exp\left\{ -{1 \over 2} \sum_{i=1}^m \left( -\beta_i+\sum_{j=1}^n \alpha_{ij} \xi_j\right)^2 \right\} \\=
&\exp\left\{-{1 \over 2}  \sum_{i=1}^m \beta_i^2 \right\} \sum_{\xi_1, \ldots, \xi_n \in \{0, 1\}} w_1^{\xi_1} \cdots w_n^{\xi_n} 
\exp\left\{ -{1 \over 2}\left( \sum_{j=1}^n \alpha_{ij} \xi_j\right)^2 \right\} \endsplit $$
and the result follows by Corollary 3.5.
{\hfill \hfill \hfill} \qed
\enddemo

Now the proof of Theorem 2.2 is obtained by rescaling.
\subhead (3.7) Proof of Theorem 2.2 \endsubhead  For $i=1, \ldots, m$ and $j=1, \ldots, n$, we define 
$$\beta_i' =\sqrt{\gamma_i}  \beta_i' \quad \text{and} \quad \alpha_{ij}' = \sqrt{\gamma_i} \alpha_{ij}.$$
Applying Corollary 3.6 to $\alpha_{ij}'$ and $\beta_i'$, we conclude that 
$$\sum_{\xi_1, \ldots, \xi_n \in \{0, 1\}} z_1^{\xi_1} \cdots z_n^{\xi_n} \exp\left\{ -{1 \over 2} \sum_{i=1}^m \gamma_i \left(-\beta_i +\sum_{j=1}^n \alpha_{ij} \xi_j \right)^2 \right\} \ne 0$$
provided 
$$\left| z_j \right| \ < \ \lambda_j \exp\left\{ -\sum_{i=1}^m \gamma_i \beta_i \alpha_{ij} \right\}=\rho_j$$
and 
$$\sqrt{\gamma_i} \sum_{j=1}^n {\lambda_j  \left| \alpha_{ij}\right|  \over 1-\lambda_j}  \ \leq \ {1 \over 2 \sqrt{c}} \quad \text{for} \quad i=1, \ldots, m.$$
{\hfill \hfill \hfill} \qed

\head 4. Concluding remarks \endhead 

\subhead (4.1) Smoothed counting of integer points \endsubhead Let $A=\left(a_{ij}\right)$ be an $m \times n$ matrix, let $b=\left(\beta_1, \ldots, \beta_m\right)$ be an 
$m$-vector and let $\gamma=\left(\gamma_1, \ldots, \gamma_m\right)$ be an $m$-vector of positive real weights. For a complex $n$-vector ${\bold z}=\left(z_1, \ldots, z_n\right)$ we introduce a series 
$$\tilde{P}_{A, b, \gamma}({\bold z})=\sum_{\xi_1, \ldots, \xi_n \in {\Bbb Z}_+} z_1^{\xi_1} \cdots z_n^{\xi_n} \exp\left\{ -\sum_{i=1}^m \gamma_i \left(-\beta_i + \sum_{j=1}^n \alpha_{ij} \xi_j \right)^2 \right\} \tag4.1.1$$
Here the external sum is taken over all $n$-tuples of non-negative integers. Clearly, (4.1.1) converges absolutely and uniformly on compact subsets of the open polydisc
$$|z_j| < 1 \quad \text{for} \quad j=1, \ldots, n,$$
although unlike (2.1.1), the function $\tilde{P}_{A, b, \gamma}$ is not a polynomial.
If we interpret $\xi_1, \ldots, \xi_n$ as independent geometric random variables such that 
$$\Pr(\xi_j =k) = (1-p_j) p_j^k \quad \text{for} \quad k=0, 1, \ldots,$$
where $0 < p_j < 1$ for $k=1, \ldots, n$, we get 
$$\EE \exp\left\{ -\sum_{i=1}^m \gamma_i \left(-\beta_i + \sum_{j=1}^n \alpha_{ij} \xi_j \right)^2 \right\}=\tilde{P}_{A, b, \gamma}\left(p_1, \ldots, p_n\right)\prod_{j=1}^n (1-p_j).$$
By more or less straightforward modification of the proof of Theorem 2.2, one can prove that $\tilde{P}_{A, b, \gamma}({\bold z}) \ne 0$ if the conditions of Theorem 2.2 are satisfied.
The proof almost repeats that of Section 3, only that in Proposition 3.1 we deal with the sum 
$$\sum_{\sigma_1, \ldots, \sigma_m =\pm 1} \prod_{j=1}^n \left(1- z_j \exp\left\{ -\sum_{i=1}^m \alpha_{ij} \sigma_i\right\}\right)^{-1}$$
and the functions $\ell(x; {\bold z})$ and $h_j(x; {\bold z})$ are replaced respectively by 
$$\tilde{\ell}(x; {\bold z}) = \prod_{j=1}^n \left(1- z_j \exp\left\{ \sqrt{-1} \sum_{i=1}^m \alpha_{ij} \sigma_i\right\}\right)^{-1}$$
and 
$$\tilde{h}_j(x; {\bold z})= {\exp\left\{ \sqrt{-1} \sum_{i=1}^m \alpha_{ij} \sigma_i\right\} \over 1- z_j \exp\left\{ \sqrt{-1} \sum_{i=1}^m \alpha_{ij} \sigma_i\right\}}.$$


\Refs
\widestnumber\key{AAAA}

\ref\key{A+99}
\by G. Ausiello, P. Crescenzi, G. Gambosi, V. Kann, A. Marchetti-Spaccamela, and M. Protasi
\book Complexity and Approximation. Combinatorial optimization problems and their approximability properties
\publ Springer-Verlag
\publaddr Berlin
\yr 1999
\endref

\ref\key{AB09}
\by S. Arora and B. Barak
\book Computational Complexity. A modern approach
\publ Cambridge University Press
\publaddr Cambridge
\yr 2009
\endref

\ref\key{AS00}
\by N. Alon and J.H. Spencer
\book The Probabilistic Method. Second edition. With an appendix on the life and work of Paul Erd\H{o}s
\bookinfo Wiley-Interscience Series in Discrete Mathematics and Optimization
\publ Wiley-Interscience [John Wiley $\&$ Sons]
\publaddr New York
\yr 2000
\endref

\ref\key{Ba16}
\by A. Barvinok
\book Combinatorics and Complexity of Partition Functions
\bookinfo  Algorithms and Combinatorics, 30
\publ Springer
\publaddr Cham
\yr 2016
\endref

\ref\key{BB21}
\by A. Barvinok and N. Barvinok
\paper More on zeros and approximation of the Ising partition function
\jour Forum of Mathematics. Sigma 
\vol 9 
\yr 2021
\pages paper no. e46, 1--18
\endref

\ref\key{BH10}
\by  A. Barvinok and J.A. Hartigan
\paper Maximum entropy Gaussian approximations for the number of integer points and volumes of polytopes
\jour  Advances in Applied Mathematics 
\vol 45 
\yr 2010
\pages no. 2, 252--289
\endref

\ref\key{BR19}
\by A. Barvinok and G. Regts
\paper Weighted counting of solutions to sparse systems of equations
\jour Combinatorics, Probability and Computing
\vol 28 
\yr 2019
\pages  no. 5, 696--719
\endref

\ref\key{B+07}
\by  M. Bayati, D. Gamarnik, D. Katz, C. Nair and P. Tetali
\paper Simple deterministic approximation algorithms for counting matchings
\inbook  STOC'07 -- Proceedings of the 39th Annual ACM Symposium on Theory of Computing
\pages 122--127
\publ ACM
\publaddr New York
\yr 2007
\endref

\ref\key{C+21}
\by Z. Chen, K. Liu and E. Vigoda
\paper Spectral independence via stability and applications to holant-type problems
\paperinfo preprint {\tt arXiv:2106.03366}
\yr 2021
\endref

\ref\key{D+14}
\by A. Dudek, M. Karpinski, A. Ruci\'nski and E. Szyma\'nska
\paper Approximate counting of matchings in $(3,3)$-hypergraphs
\bookinfo Algorithm theory -- SWAT 2014
\pages  380--391
\inbook Lecture Notes in Computer Science 
\vol 8503
\publ Springer
\publaddr Cham
\yr 2014
\endref

\ref\key{FV18}
\by S. Friedli and Y. Velenik
\book Statistical Mechanics of Lattice Systems. A concrete mathematical introduction
\publ Cambridge University Press
\publaddr Cambridge
\yr  2018
\endref

\ref\key{Ga20}
\by D. Gamarnik
\paper Correlation decay and the absence of zeros property of partition functions
\paperinfo preprint {\tt arXiv:2011.04915}
\yr 2020
\endref

\ref\key{JS89}
\by M. Jerrum and A, Sinclair
\paper Approximating the permanent
\jour SIAM Journal on Computing
\vol 18 
\yr 1989
\pages no. 6, 1149--1178
\endref

\ref\key{J+04}
\by M. Jerrum, A. Sinclair and E. Vigoda
\paper A polynomial-time approximation algorithm for the permanent of a matrix with nonnegative entries
\jour Journal of the ACM
\vol  51
\yr 2004
\pages no. 4, 671--697
\endref

\ref\key{Ke18}
\by P. Keevash
\paper Counting designs
\jour Journal of the European Mathematical Society (JEMS) 
\vol 20 
\yr 2018
\pages  no. 4, 903--927
\endref

\ref\key{Kr01}
\by S.G. Krantz
\book Function Theory of Several Complex Variables. Reprint of the 1992 edition
\publ AMS Chelsea Publishing
\publaddr Providence, RI
\yr 2001
\endref

\ref\key{LL13}
\by N. Linial and Z. Luria
\paper An upper bound on the number of Steiner triple systems
\jour Random Structures $\&$ Algorithms 
\vol 43 
\yr 2013
\pages 4, 399--406
\endref

\ref\key{LL14}
\by N. Linial and Z. Luria
\paper An upper bound on the number of high-dimensional permutations
\jour Combinatorica 
\vol 34 
\yr 2014
\pages no. 4, 471--486
\endref

\ref\key{LY52}
\by T.D. Lee and C.N. Yang
\paper Statistical theory of equations of state and phase transitions. II. Lattice gas and Ising model
\jour  Physical Review (2) 
\vol 87 
\yr 1952
\pages  410--419
\endref

\ref\key{MR95}
\by R. Motwani and P. Raghavan
\book Randomized Algorithms
\publ Cambridge University Press
\publaddr Cambridge
\yr 1995
\endref

\ref\key{NN94}
\by Yu. Nesterov and A. Nemirovskii
\book Interior-Point Polynomial Algorithms in Convex Programming
\bookinfo SIAM Studies in Applied Mathematics, {\bf 13}
\publ Society for Industrial and Applied Mathematics (SIAM)
\publaddr Philadelphia, PA
\yr 1994
\endref

\ref\key{PR17}
\by V. Patel and G. Regts
\paper Deterministic polynomial-time approximation algorithms for partition functions and graph polynomials
\jour SIAM Journal on Computing 
\vol 46 
\yr 2017
\pages no. 6, 1893--1919
\endref

\ref\key{PP20}
\by I. Pak and G. Panova
\paper Bounds on Kronecker coefficients via contingency tables
\jour Linear Algebra and its Applications 
\vol 602 
\yr 2020
\pages 157--178
\endref

\ref\key{Po18}
\by V.N. Potapov
\paper On the number of SQSs, Latin hypercubes and MDS codes
\jour  Journal of Combinatorial Designs
\vol 26 
\yr 2018
\pages  no. 5, 237--248
\endref

\ref\key{PR17}
\by V. Patel and G. Regts
\paper Deterministic polynomial-time approximation algorithms for partition functions and graph polynomials
\jour SIAM Journal on Computing 
\vol 46 
\yr 2017
\pages no. 6, 1893--1919
\endref 

\ref\key{S+19}
\by R. Song, Y. Yin and J. Zhao
\paper Counting hypergraph matchings up to uniqueness threshold
\jour Information and Computation 
\vol 266 
\yr 2019
\pages  75--96
\endref

\ref\key{Va79}
\by L.G. Valiant
\paper The complexity of enumeration and reliability problems
\jour SIAM Journal on Computing 
\vol 8 
\yr 1979
\pages no. 3, 410--421
\endref 


\endRefs
\enddocument
\end